\documentclass[12pt,reqno]{amsart}
\usepackage{amssymb,amsthm}
\usepackage{braket,caption,dsfont,enumerate,esint,float,mathdots,mathtools}
\usepackage{perpage,pdfpages,pgfplots,thmtools}
\usepackage[margin=2.2cm]{geometry}

\usepackage[T1]{fontenc}
\usepackage[english]{babel}
\usepackage{hyperref,lmodern,microtype}

\allowdisplaybreaks 
\pgfplotsset{compat=1.17}
\MakePerPage{footnote}

\DeclareMathOperator{\Area}{Area}

\DeclareMathOperator{\dist}{dist}

\DeclareMathOperator{\Op}{Op}

\DeclareMathOperator{\supp}{supp}

\newcommand{\C}{\mathbb{C}}

\newcommand{\comp}{\mathrm{comp}}

\newcommand{\E}{\mathbb{E}}

\newcommand{\Gc}{\mathcal{G}}

\newcommand{\Hb}{\mathbb{H}}

\newcommand{\indicator}{\mathds{1}}

\newcommand{\M}{\mathbb{M}}

\newcommand{\N}{\mathbb{N}}
\newcommand{\Nc}{\mathcal{N}}

\newcommand{\ol}{\overline}

\newcommand{\R}{\mathbb{R}}

\newcommand{\Sb}{\mathbb{S}}

\newcommand{\SH}{\mathbb{SH}}

\newcommand{\SO}{\mathrm{SO}}

\newcommand{\SU}{\mathrm{SU}}

\newcommand{\T}{\mathbb{T}}

\newcommand{\Z}{\mathbb{Z}}

\theoremstyle{plain}
\newtheorem{thm}{Theorem}[section]

\newtheorem{cor}[thm]{Corollary}
\newtheorem{lemma}[thm]{Lemma}
\newtheorem{pblm}[thm]{Problem}
\newtheorem{prop}[thm]{Proposition}

\theoremstyle{definition}
\newtheorem*{defn}{Definition}
\newtheorem*{ex}{Example}

\theoremstyle{remark}
\newtheorem*{rmk}{Remark}

\numberwithin{equation}{section}

\title{Spherical harmonics and point configurations on the sphere}
\author{Xiaolong Han}
\email{xiaolong.han@csun.edu}
\address{Department of Mathematics, California State University, Northridge, CA 91330, USA}



\dedicatory{Dedicated to the memory of Steve Zelditch}
\thanks{} 

\begin{document}
\maketitle

\begin{abstract}
We construct spherical harmonics on the two-dimensional unit sphere as superpositions of Gaussian beams whose poles form well-separated point configurations. The distributional and analytic properties of the resulting spherical harmonics are determined by the geometry of these poles: when the configuration is equidistributed, the sequence of harmonics exhibits quantum ergodicity, while their $L^\infty$ norms are quantitatively controlled by the maximal clustering of poles within small neighborhoods of great circles.
\end{abstract}

\section{Introduction}
A problem of both theoretical and practical importance is that of distributing a large number of points on the unit sphere $\Sb\subset\R^3$ in a well-balanced manner; see Hardin-Michaels-Saff \cite{HMS} for an overview of classical and contemporary point configurations arising from diverse applications. In this paper, we establish a new connection between point configurations on $\Sb$ and the construction of spherical harmonics exhibiting prescribed (and often extreme) geometric and analytic properties.

\textit{Spherical harmonics} on $\Sb$ are the restrictions to $\Sb$ of homogeneous harmonic polynomials in $\R^3$. They arise naturally in mathematics and physics as the eigenfunctions of the Laplacian $\Delta_\Sb$ on $\Sb$. For $N\in\N$, denote by $\SH_N$ the space of spherical harmonics of degree $N$. Then 
$$-\Delta_\Sb u=N(N+1)u\quad\text{for }u\in\SH_N.$$
The study of spherical harmonics thus falls within the broader theory of Laplacian eigenfunctions on compact manifolds. We focus on their distribution and $L^p$ norms, particularly in the high-energy limit as eigenvalues tend to infinity. These behaviors reflect the underlying geometry and the dynamics of the geodesic flow, and have been extensively investigated; see Sogge \cite{So} and Zelditch \cite{Ze3, Ze4} for comprehensive surveys.

On $\Sb$, the eigenspace $\SH_N$ has dimension $\dim\SH_N=2N+1$ and contains a wide range of eigenfunctions displaying distinct geometric and analytic behaviors. This large multiplicity makes spherical harmonics a central model for studying localization and delocalization phenomena, as well as for testing the sharpness of $L^p$ estimates in the general theory of Laplacian eigenfunctions. An important example is provided by the Gaussian beams.

\begin{ex}[Gaussian beams]
Let $x=(x_1,x_2,x_3)\in\Sb$. Define
\begin{equation}\label{eq:Q0}
Q_0(x)=C_NN^\frac14\left(x_1+ix_2\right)^N,
\end{equation}
where $C_N>0$ is chosen so that $\|Q_0\|_{L^2(\Sb)}=1$. There exist absolute constants $c_1,c_2>0$ such that $c_1\le C_N\le c_2$ for all $N\in\N$; see Zelditch \cite[Section 4.4.5]{Ze4}. Note that $Q_0$ depends on $N$, and there is one such spherical harmonic in $\SH_N$ for each $N\in\N$. For simplicity of notation, we omit the explicit dependence on $N$.

Write the spherical coordinates
$$x=\left(x_1,x_2,x_3\right)=(\sin\phi\cos\theta,\sin\phi\sin\theta,\cos\phi),$$
where $\phi\in[0,\pi]$ and $\theta\in[0,2\pi)$. Then
\begin{equation}\label{eq:Q0polar}
Q_0(\phi,\theta)=C_NN^\frac14(\sin\phi)^Ne^{iN\theta}=C_NN^\frac14(\cos\alpha)^Ne^{iN\theta},
\end{equation}
where $\alpha=\frac\pi2-\phi$ is the angle between $x$ and the equator $G_0=\{\phi=\frac\pi2\}$. It follows that $Q_0$ is localized in an $O(N^{-\frac12})$-neighborhood of $G_0$ and decays exponentially in the transverse direction. Moreover, $Q_0$ oscillates at frequency $N$ in the counterclockwise direction by the right-hand rule with respect to the north pole $p_0$ (i.e., $\phi=0$), which we call the \textit{pole} of $Q_0$. Since $Q_0$ oscillates along $G_0$ and has Gaussian decay in the transverse direction, it is commonly referred to as a \textit{Gaussian beam} (or a \textit{highest weight spherical harmonic}).

Any rotation of $Q_0$ by an element of $\SO(3)$ remains a spherical harmonic in $\SH_N$ (and we choose the rotation to be orientation-preserving). Thus, for any $p\in\Sb$, there exists a Gaussian beam with pole $p$. In fact, there is a family of Gaussian beams sharing the same pole $p$; these functions differ only by a phase shift. See Han \cite[Section 2]{Ha1} for a detailed discussion. In particular, the Gaussian beam
$$C_NN^\frac14\left(x_1-ix_2\right)^N=C_NN^\frac14(\sin\phi)^Ne^{-iN\theta}$$
has the south pole (i.e., $\phi=\pi$) as its pole and is orthogonal to $Q_0$. It exhibits the same localization around the equator $G_0$ as $Q_0$ but oscillates in the opposite direction along $G_0$.
\end{ex}

Gaussian beams serve as the basic building blocks in our construction of spherical harmonics. For a point configuration $\{p_j\}_{j=1}^m\subset\Sb$, define $Q_j$ as the Gaussian beam with pole $p_j$, $j=1,...,m$. Set
\begin{equation}\label{eq:uN}
F_N=\sum_{j=1}^mQ_j\quad\text{and}\quad u_N=\frac{F_N}{\left\|F_N\right\|_{L^2(\Sb)}}.
\end{equation}
Throughout, we assume that $m\le\dim\SH_N=2N+1$, and let $m=m(N)\to\infty$. When the point configurations $\{p_j\}_{j=1}^m$ satisfy the following structural conditions, the corresponding spherical harmonics $u_N$ in \eqref{eq:uN} exhibit the desired distributional or $L^\infty$ properties in Theorems \ref{thm:QE} and \ref{thm:BddSHs} below.

\begin{defn}[Conditions on point configurations]
Let $\{p_j\}_{j=1}^m\subset\Sb$, $m\in\N$, be a family of point configurations, one for each $m$; the dependence of the points $p_j$ on $m$ is suppressed in the notation. Conditions (S), (E), and (G) below are imposed on the family as $m\to\infty$.

(S). [\textit{Separation}] There exists an absolute constant $c>0$ such that
\begin{equation}\label{eq:S}
\dist\left(p_j,p_k\right)\ge\frac{c}{\sqrt m}\quad\text{for all }j\ne k.
\end{equation}

(E). [\textit{Equidistribution}] For any $f\in C^\infty(\Sb)$,
\begin{equation}\label{eq:E}
\lim_{m\to\infty}\frac1m\sum_{j=1}^m f\left(p_j\right)=\frac{1}{4\pi}\int_\Sb f(x)\,dx.
\end{equation}

(G). [\textit{Bounded clustering around great circles}] There is a function $L(m)\ge1$ such that
\begin{equation}\label{eq:G}
\#\left\{j:p_j\in\Nc_{m^{-1}}(G)\right\}\le L(m)\quad\text{for all great circles }G\subset\Sb.
\end{equation}

Here $\dist$ denotes the geodesic distance on $\Sb$, and $\Nc_r(\Omega)=\{x\in\Sb:\dist(\Omega,x)\le r\}$ denotes the $r$-neighborhood of $\Omega\subset\Sb$.
\end{defn}

\begin{rmk}
We note that for any point configuration $\{p_j\}_{j=1}^m$, the scale $m^{-\frac12}$ in \eqref{eq:S} is optimal: it is impossible to place $m$ points on $\Sb$ that remain separated by a distance of at least $m^{-s}$ for some $s<\frac12$.
\end{rmk}

\begin{thm}[Quantum ergodic spherical harmonics]\label{thm:QE}
Suppose that $m=O(N^{\rho_0})$ for some $0<\rho_0<1$. If (S) and (E) hold, then the spherical harmonics $u_N$ in \eqref{eq:uN} are quantum ergodic as $N\to\infty$.
\end{thm}

\begin{thm}[Spherical harmonics with controlled $L^\infty$ bounds]\label{thm:BddSHs}
Suppose that $m=\lfloor\delta N^\frac12\rfloor+1$ for a sufficiently small $\delta>0$. If (S) and (G) hold, then the spherical harmonics $u_N$ in \eqref{eq:uN} satisfy
$$\|u_N\|_{L^\infty(\Sb)}\le C\delta^{-\frac12}L(m),$$
where $C>0$ is an absolute constant. In particular, if (G) holds with $L(m)=O(1)$ as $m\to\infty$, then $u_N$ in \eqref{eq:uN} are uniformly bounded as $N\to\infty$.
\end{thm}

\begin{rmk}
The relationship between point configurations and the construction of uniformly bounded spherical harmonics in Theorem \ref{thm:BddSHs} was observed by Bourgain as early as 2013; his question is recorded in \cite[Page 957]{Bo3} and remains unpublished. We arrived at it independently. We are grateful to C. Demeter and P. Varj\'u for drawing our attention to this observation.
\end{rmk}

Next, we present concrete examples of point configurations that satisfy the above conditions. We then place Theorems \ref{thm:QE} and \ref{thm:BddSHs} in the broader context of Laplacian eigenfunction theory.

\subsection{Point configurations on the sphere}\label{sec:pt}
Constructing large collections of points that are ``well distributed'' according to various criteria has numerous applications. For a comprehensive overview of this topic, see Saff-Kuijlaars \cite{SK}. 

Quantum ergodic spherical harmonics in Theorem \ref{thm:QE} are constructed from point configurations satisfying the (S) separation and (E) equidistribution conditions. Several examples of such configurations are known; we describe two representative cases below.

\begin{ex}[Spherical designs]
Let $t\in\N$. A point configuration $\{p_j\}_{j=1}^m$ is called a \textit{$t$-spherical design} if
$$\frac1m\sum_{j=1}^mP\left(p_j\right)=\int_\Sb P(x)\,dx$$
for all polynomials $P$ of degree at most $t$. Bondarenko-Radchenko-Viazovska \cite{BRV1} proved the existence of $t$-spherical designs $\{p_j\}_{j=1}^m$ with $m\le ct^2$, where $c>0$ is an absolute constant. Since any $f\in C^\infty(\Sb)$ can be approximated by polynomials, these configurations automatically satisfy (E) equidistribution. Moreover, \cite{BRV2} established the existence of $t$-spherical designs that also satisfy (S) separation.

Their results \cite{BRV1, BRV2} extend to higher-dimensional spheres; see \cite{BRV1, BRV2} for further background and related developments on spherical designs.
\end{ex}

The argument for the $t$-spherical designs in \cite{BRV1, BRV2} is non-constructive, and producing an explicit $t$-spherical design for large $t$ remains an open problem. We now describe an explicit family of point configurations satisfying both (S) separation and (E) equidistribution.

\begin{ex}[Point configurations from Hecke operators]
Lubotzky-Phillips-Sarnak \cite{LPS1, LPS2} constructed large families of point configurations on $\Sb$. We briefly recall one approach using Hecke operators. Let $\Hb$ denote the Hamiltonian quaternions. For $\alpha=\alpha_0+\alpha_1 i+\alpha_2j+\alpha_3k\in\Hb(\R)$, define
$$S_\alpha=\frac{1}{|\alpha|}\begin{pmatrix}
\alpha_0+\alpha_1i & \alpha_2+\alpha_3i\\
-\alpha_2+\alpha_3i & \alpha_0-\alpha_1i\end{pmatrix}\in\SU(2),$$
where
$$|\alpha|^2=\alpha\ol\alpha=\alpha_0^2+\alpha_1^2+\alpha_2^2+\alpha_3^2.$$
The matrix $S_\alpha$ acts on $\C\cup\{\infty\}$ by a linear fractional transformation, which via stereographic projection defines a rotation in $\SO(3)$ on $\Sb$. They showed that the orbit
$$\left\{S_\alpha(x):\alpha\in\Hb(\Z),|\alpha|^2=n\right\}$$
satisfies (S) separation and (E) equidistribution under suitable conditions on $n$. The proof relies on the spectral analysis of the Hecke operator
\begin{equation}\label{eq:Hecke}
T_nf=\sum_{\alpha\in\Hb(\Z),|\alpha|^2=n}f\circ S_\alpha.
\end{equation}
For example, when $n=5$, $T_5$ corresponds to six rotations through an angle $\arccos(-\frac35)$ about the $x_1$, $x_2$, and $x_3$ axes. The absolute value of the second-largest eigenvalue of $T_5$ is bounded by $2\sqrt5$, which implies strong distribution properties for the orbits generated by $\{S_\alpha:\alpha\in\Hb(\Z),|\alpha|^2=5^k\}$ as $k\to\infty$; see \cite{LPS1, LPS2} for details.
\end{ex}

\begin{figure}[H]
\begin{tikzpicture}
\draw[thick] (0,0) circle (4);
\draw[thick] (-4,0) arc (180:360:4 and 0.6);
\draw[thick,dashed] (4,0) arc (0:180:4 and 0.6);
\draw[] (-3.967,0.5) arc (180:360:3.967 and 0.6);
\draw[dashed] (3.967,0.5) arc (0:180:3.967 and 0.6);
\draw[] (-3.967,-0.5) arc (180:360:3.967 and 0.6);
\draw[dashed] (3.967,-0.5) arc (0:180:3.967 and 0.6);
\draw (4,0) node[right]{$G$};
\draw[thick,<->] (-4.2,-0.5) arc (190:175:4.2) node[left,midway]{$2m^{-1}$};
\foreach \i in {1,...,500} {
  \pgfmathparse{360*rnd}
  \let\angle\pgfmathresult
  \pgfmathparse{4*sqrt(rnd)}
  \let\radius\pgfmathresult
  \draw[blue] (\angle:\radius) node[]{$*$};
}
\end{tikzpicture}
\caption{Condition (G) bounded clustering around great circles}
\label{fig:pt}
\end{figure}

Condition (G) bounded clustering around great circles \eqref{eq:G} is of a different nature, in that it constrains the configuration at the scale $m^{-1}$ relative to every great circle; see Figure \ref{fig:pt}.

Observe that, by the rotational symmetry of $\Sb$, the $m^{-1}$-neighborhood $\Nc_{m^{-1}}(G)$ of any great circle $G$ has the same area as $\Nc_{m^{-1}}(G_0)=\{|\alpha|\le m^{-1}\}$, that is,
$$\Area\left(\Nc_{m^{-1}}(G)\right)=\int_0^{2\pi}\int_{-\frac1m}^{\frac1m}\cos\alpha\,d\alpha d\theta=4\pi\sin\frac1m\le\frac{4\pi}{m}.$$
Hence, averaged over all great circles $G$, the number of points from $\{p_j\}_{j=1}^m$ contained in $\Nc_{m^{-1}}(G)$ is uniformly bounded. Therefore, \eqref{eq:G} with $L(m)=O(1)$ is natural in a certain averaged sense. Constructing an explicit configuration for which the number of points contained in $\Nc_{m^{-1}}(G)$ remains uniformly bounded for \emph{all} great circles is, however, a subtle and challenging problem. 

While we do not provide such a configuration here, in Appendix \ref{sec:prob} we construct point configurations satisfying (S), (E), and (G) with $L(m)=O\left(\frac{\log m}{\log\log m}\right)$. This is achieved by introducing independent random perturbations to configurations already satisfying (S) and (E), such as those arising from Hecke operators. We also refer to Demeter-Zhang \cite{DZ} for a recent discussion in a planar model, namely, large collections of points in a square avoiding clustering around all straight lines, though the spherical case considered here may be more subtle.

\subsection{Quantum ergodic spherical harmonics}\label{sec:QE0}
Quantum ergodic eigenfunctions are equidistributed in the phase space of the manifold. The study of such distributional properties originated on manifolds whose geodesic flow is ergodic, notably on compact hyperbolic manifolds, that is, compact manifolds with constant curvature $-1$; see \cite{KH}. These properties are naturally described using semiclassical measures, which we briefly recall.

Let $\M$ be a compact manifold and $T^*\M=\{(x,\xi):x\in\M,\xi\in T_x^*\M\}$ its cotangent bundle. Fix $h_0\in(0,1)$ and let $h\in(0,h_0)$ denote the semiclassical parameter. Consider eigenfunctions of the semiclassical Laplacian
$$\left(h^2\Delta_\M+1\right)u_h=0.$$
A sequence of eigenfunctions $\{u_{h_k}\}_{k=1}^\infty\subset L^2(\M)$ is said to induce a \textit{semiclassical measure} $\mu$ on $T^*\M$ if, for every $a\in C_0^\infty(T^*\M)$,
$$\langle\Op_{h_k}(a)u_{h_k},u_{h_k}\rangle\to\int_{T^*\M}a\,d\mu\quad\text{as }h_k\to0,$$
where $\Op_h(a)$ denotes the semiclassical pseudo-differential operator with symbol $a$. If $a=a(x)\in C^\infty(\M)$, then $\Op_h(a)=a$, and 
$$\left\langle\Op_h(a)u_h,u_h\right\rangle=\int_\M a\left|u_h\right|^2.$$ 
Thus, the semiclassical measure can be viewed as the lift of the measure $|u_h|^2$ from the physical space $\M$ to the phase space $T^*\M$, capturing the joint distribution of $u_h$ in position and frequency.

It is well known that any semiclassical measure is supported on the cosphere bundle
$$S^*\M=\left\{(x,\xi)\in T^*\M:|\xi|_x=1\right\},$$
and is invariant under the geodesic flow on $S^*\M$; see Zworski \cite[Section~5.2]{Zw}. Among all invariant probability measures, the Liouville measure $\mu_L$, normalized so that $\mu_L(S^*\M)=1$, provides the canonical uniform measure on $S^*\M$.

The celebrated Quantum Ergodicity Theorem of Shnirelman \cite{Shn}, Zelditch \cite{Ze1}, and Colin de Verdi\`ere \cite{CdV} asserts that if the geodesic flow on $S^*\M$ is ergodic, then every eigenbasis contains a full-density subsequence whose semiclassical measure equals $\mu_L$. Having $\mu_L$ as the limit of their semiclassical measures has led to such eigenfunctions being called \textit{equidistributed in the phase space}; they are also commonly said to be \textit{quantum ergodic}, since Laplacian eigenfunctions represent stationary states of the quantum system associated with the geodesic flow.

We note that quantum ergodicity implies equidistribution on $\M$, but is strictly stronger than this property.
\begin{ex}[Toral eigenfunctions]
On the flat torus $\T^d=\R^d/2\pi\Z^d$, consider the eigenfunctions
$$u_{h_k}(x)=e^{i\langle\lambda_k,x\rangle}=e^{\frac{i}{h_k}\left\langle\frac{\lambda_k}{\left|\lambda_k\right|},x\right\rangle},$$
with semiclassical parameter $h_k=|\lambda_k|^{-1}$. If $\frac{\lambda_k}{|\lambda_k|}\to\xi_0$ as $|\lambda_k|\to\infty$ for some fixed $\xi_0\in\Sb^{d-1}$, then 
$$\left\langle\Op_{h_k}(a)u_{h_k},u_{h_k}\right\rangle\to\int_{\T^d}a\left(x,\xi_0\right)\,dx\quad\text{as }h_k\to0.$$
The corresponding semiclassical measure is $dx\delta_{\xi=\xi_0}$, supported on the invariant set $\T^d\times\{\xi_0\}\subset S^*\T^d$. This illustrates that although the toral eigenfunctions $e^{i\langle\lambda_k,x\rangle}$ are equidistributed in the physical space $\T^d$, they are completely localized in frequency and therefore not quantum ergodic.
\end{ex}

\begin{rmk}[Quantum Unique Ergodicity]
Rudnick-Sarnak’s Quantum Unique Ergodicity (QUE) conjecture \cite{RS} asserts that the entire eigenbasis on a compact hyperbolic manifold is quantum ergodic. On arithmetic hyperbolic surfaces, one can choose a joint eigenbasis of the Laplacian and the Hecke operators, reflecting the arithmetic symmetries of the surface. The QUE conjecture has been established for this Hecke eigenbasis by Lindenstrauss \cite{L} and Brooks-Lindenstrauss \cite{BrLi}. 
\end{rmk}

\begin{rmk}[Hecke spherical harmonics]
Recall the Hecke operators \eqref{eq:Hecke} on $\Sb$, which generate point configurations satisfying (S) separation and (E) equidistribution. These configurations are then used to construct quantum ergodic spherical harmonics in Theorem \ref{thm:QE}.

The Hecke operators can also directly define joint eigenfunctions with the Laplacian $\Delta_\Sb$, which are called Hecke spherical harmonics. An analogue of the QUE conjecture for Hecke spherical harmonics on $\Sb$ was proposed by B\"ocherer-Sarnak-Schulze-Pillot \cite[Conjecture 1]{BSSP} and remains open. The precise relation between the quantum ergodic spherical harmonics in Theorem \ref{thm:QE} arising from ``Hecke points'' and the Hecke spherical harmonics is also unclear, and we leave it for future study.
\end{rmk}

In the semiclassical framework, spherical harmonics $u\in\SH_N$ satisfy
$$\left(h_N^2\Delta_{\Sb}+1\right)u=0,$$ 
with semiclassical parameter
$$h_N=\frac{1}{\sqrt{N(N+1)}}\approx N^{-1}\to0\quad\text{as }N\to\infty.$$ 
The geodesic flow on $\Sb$ is completely integrable and not ergodic \cite{KH}. As a result, there exist spherical harmonics that exhibit strong localization and therefore fail to be quantum ergodic.

\begin{ex}[Semiclassical measure of Gaussian beams]
Let $Q_p$ be a Gaussian beam with a fixed pole $p\in\Sb$. Then
$$\lim_{N\to\infty}\left\langle\Op_h(a)Q_p,Q_p\right\rangle=\frac{1}{2\pi}\int_{G_p}a\left(x,\xi_p\right)\,dl,$$
where $\xi_p\in S^*G_p$ is the covector determined by the right-hand rule with respect to $p$ and $dl$ denotes the arclength measure on $G_p$. Hence, the corresponding semiclassical measure is $\frac{1}{2\pi}dl\delta_{\xi=\xi_p}$, supported on the invariant set $G_p\times\{\xi_p\}\subset S^*\Sb$. See Proposition \ref{prop:smG} for details.
\end{ex}

Nevertheless, there exist spherical harmonics that are quantum ergodic despite the integrable nature of the geodesic flow on $\Sb$. Jakobson-Zelditch \cite{JZ} showed that every invariant probability measure on $S^*\Sb$ arises as the semiclassical measure of some sequence of spherical harmonics. In particular, this includes the Liouville measure $\mu_L$, establishing the existence of quantum ergodic spherical harmonics. Moreover, Brooks-Le Masson-Lindenstrauss \cite{BLML} proved that a joint eigenbasis of the Laplacian and a free group generated by finitely many rotations in $\SO(3)$ contains a full-density subsequence of quantum ergodic spherical harmonics.

Since $\dim\SH_N=2N+1$, one can introduce a natural probability measure on the space of eigenbases. In this probabilistic framework, a random eigenbasis of $\SH_N$ satisfies quantum unique ergodicity almost surely; see Zelditch \cite{Ze2}, VanderKam \cite{V}, and Burq-Lebeau \cite{BuLe} for precise formulations.

Our result in Theorem \ref{thm:QE} complements existing existence and probabilistic results by providing explicit constructions of quantum ergodic spherical harmonics based on point configurations satisfying the structural conditions of (S) separation and (E) equidistribution. There are many examples of point configurations $\{p_j\}_{j=1}^m$ satisfying (S) and (E); see Section \ref{sec:pt}. Moreover, Theorem \ref{thm:QE} remains valid as long as $m = O(N^{\rho_0})$ for any $\rho_0 \in (0,1)$. Therefore, our result yields a large and explicit family of quantum ergodic spherical harmonics.

\subsection{Spherical harmonics with controlled $L^\infty$ bounds}
Spherical harmonics satisfy the $L^\infty$ bound for eigenfunctions on compact manifolds established by H\"ormander \cite{Ho}:
\begin{equation}\label{eq:Linfty}
\|u\|_{L^\infty(\Sb)}\le CN^{\frac12}\|u\|_{L^2(\Sb)}\quad\text{for all }u\in\SH_N,
\end{equation} 
where $C>0$ is an absolute constant.

A fundamental and largely open problem is to identify and characterize eigenfunctions with minimal $L^\infty$ growth, namely those satisfying $\|u\|_{L^\infty}\le C\|u\|_{L^2}$; that is, ($L^2$-normalized) eigenfunctions that remain uniformly bounded as the eigenvalue tends to infinity. (The reverse inequality $\|u\|_{L^2}\le C\|u\|_{L^\infty}$ follows trivially from H\"older’s inequality on compact manifolds.)

\begin{ex}[Toral eigenfunctions]
Consider the flat tori $\T^d=\R^d/\Lambda$, where $\Lambda$ is a lattice, for example $2\pi\Z^d$ as above. There exists a basis of uniformly bounded eigenfunctions of the form $e^{i\langle\lambda,x\rangle}$, where $x\in\T^d$ and $\lambda\in\Lambda^\star$, the dual lattice of $\Lambda$. One may also consider eigenfunctions on a fundamental domain of $\R^d/\Lambda$ with appropriate boundary conditions.
\end{ex}

Toth-Zelditch \cite{TZ} showed that, under suitable conditions including ``quantum integrability,'' a compact Riemannian manifold in which every eigenbasis is uniformly bounded must be flat. See \cite{TZ} for further results and a comprehensive discussion.

Aside from flat tori, the only manifolds known to admit uniformly bounded eigenfunctions are certain higher-dimensional spheres. In particular, $\Sb^3$ can be realized as the boundary of the unit ball in the complex space $\C^2$. Bourgain constructed an orthonormal basis of uniformly bounded homogeneous holomorphic polynomials, which are also spherical harmonics on $\Sb^3$ \cite{Bo1}. He subsequently extended this construction to $\Sb^5$, regarded as the boundary of the unit ball in $\C^3$ \cite{Bo2}. On more general compact K\"ahler manifolds, Shiffman \cite{Shi} and Marzo-Ortega-Cerd\'a \cite{MOC} proved the existence of uniformly bounded holomorphic sections, although these do not form a full basis in the complex domain. Applied to unit spheres in complex space, these results imply the existence of uniformly bounded spherical harmonics of arbitrary degree on odd-dimensional spheres. However, the methods in \cite{Bo1, Shi, Bo2, MOC} do not extend to even-dimensional spheres.

Our result in Theorem \ref{thm:BddSHs} introduces a novel approach for constructing spherical harmonics on $\Sb$, in which the $L^\infty$ bound is directly controlled by the quantity $L(m)$ in \eqref{eq:G}, representing the maximal number of points contained in the $m^{-1}$-neighborhoods of great circles among the $m$ points used in the construction. Consequently, constructing uniformly bounded spherical harmonics on $\Sb$ reduces to finding well-separated point configurations with bounded clustering around great circles.

\begin{rmk}
The $L^\infty$ estimates for $u_N$ in Theorem \ref{thm:BddSHs} rely crucially on the geometric assumptions (S) and (G) on the poles. By arranging the poles in different configurations, one can produce a wide variety of behaviors for the resulting spherical harmonics. For instance, placing roughly $N^{\frac12}$ poles evenly along the equator with spacing $\approx N^{-\frac12}$ yields a spherical harmonic of the form \eqref{eq:uN} that attains the maximal $L^\infty$ growth of order $N^{\frac12}$ in \eqref{eq:Linfty}; see Guo-Han-Tacy \cite[Section 4.2]{GHT}.
\end{rmk}

We note that for each $N\in\N$, Theorem \ref{thm:BddSHs} produces a single spherical harmonic in $\SH_N$, which therefore has density zero as $N\to\infty$ within an eigenbasis, since $\dim\SH_N=2N+1$. The spherical harmonics constructed in \cite{Bo1, Shi, Bo2, MOC} using complex-analytic methods likewise have density zero within an eigenbasis, as the dimension of the holomorphic functions grows more slowly than that of the full space of spherical harmonics. For instance, the space of spherical harmonics of degree $N$ on $\Sb^3$ has dimension of order $N^2$, while the holomorphic polynomials of degree $N$ on $\C^2$ have dimension $N$. Consequently, the uniformly bounded spherical harmonics in \cite{Bo1} have density $O(N^{-1})\to0$ as $N\to\infty$ within an eigenbasis. This naturally leads to the following question, which is completely open in all dimensions.

\begin{pblm}[Uniformly bounded eigenbasis of spherical harmonics]
Determine whether there exist uniformly bounded spherical harmonics with asymptotically positive density in an eigenbasis, and whether an entire eigenbasis can consist of uniformly bounded spherical harmonics.
\end{pblm}

\subsection{Quantum ergodic and uniformly bounded spherical harmonics}
The constructions underlying Theorems \ref{thm:QE} and \ref{thm:BddSHs} can be combined to yield the following corollary.

\begin{cor}[Quantum ergodic and uniformly bounded spherical harmonics]\label{cor}
Suppose that $m=\lfloor\delta N^\frac12\rfloor+1$ for a sufficiently small $\delta>0$. If the conditions (S), (E), and (G) hold and $L(m)=O(1)$ as $m\to\infty$, then the spherical harmonics $u_N$ in \eqref{eq:uN} are both quantum ergodic and uniformly bounded as $N\to\infty$.
\end{cor}

We place these special spherical harmonics in the broader context of eigenfunction behaviors to highlight their distinctiveness.

\begin{ex}[Toral eigenfunctions]
On the tori $\T^d=\R^d/2\pi\Z^d$, the eigenfunctions $e^{ik\cdot x}$ are uniformly bounded. However, as discussed earlier, they are \emph{not} quantum ergodic. Moreover, other uniformly bounded eigenfunctions, such as $\cos(Nx_1)\sin(x_2)$ as $N\to\infty$ on $\T^2$, also fail quantum ergodicity; see Jakobson \cite{J}.
\end{ex}

\begin{ex}[Hecke eigenfunctions]
On arithmetic hyperbolic surfaces, the Hecke eigenfunctions are known to be quantum ergodic \cite{BrLi, L}. Nevertheless, Iwaniec-Sarnak \cite[Theorem 0.1(b)]{IS} proved that there is an absolute constant $c>0$ such that $\|u\|_{L^\infty}\ge c\sqrt{\log\log\lambda}$ along a subsequence of Hecke eigenfunctions $u$ (normalized in $L^2$), where $\lambda\to\infty$ denotes the eigenvalue of $u$. Hence, the Hecke eigenfunctions are not uniformly bounded.
\end{ex}

\begin{ex}[Uniformly bounded spherical harmonics on $\Sb^3$]
Bourgain \cite{Bo1} constructed uniformly bounded spherical harmonics on $\Sb^3$. The author recently showed that these are not quantum ergodic; indeed, their semiclassical measures are supported on the family of Clifford tori \cite{Ha2}. It remains an open question, but it is plausible that the uniformly bounded spherical harmonics on $\Sb^5$ constructed in \cite{Bo2} exhibit a similar localization phenomenon.
\end{ex}

\begin{ex}[Random spherical harmonics]
A random eigenbasis of spherical harmonics is known to satisfy quantum unique ergodicity almost surely \cite{Ze2, V, BuLe}. On the other hand, $\|u\|_{L^\infty(\Sb)}\approx\sqrt{\log N}$ almost surely \cite{V, BuLe}.
\end{ex}

The spherical harmonics in Corollary \ref{cor} provide the first examples of Laplacian eigenfunctions that are simultaneously quantum ergodic and uniformly bounded. We stress, however, that this result is \textit{conditional} on the construction of specific point configurations.

In Appendix \ref{sec:prob}, we construct point configurations satisfying (S), (E), and (G) with $L(m)$ of logarithmic order in $m$. Applying the same argument as in Corollary \ref{cor} then yields the following \textit{unconditional} result.

\begin{thm}\label{thm:log}
There are quantum ergodic spherical harmonics $u_N$, $u_N\in\SH_N$, such that
$$\left\|u_N\right\|_{L^\infty(\Sb)}\le\frac{C\log N}{\log\log N},$$
where $C>0$ is an absolute constant.
\end{thm}

\subsection*{Acknowledgements}
I would like to thank R. Zhang for several stimulating discussions on Bourgain’s work \cite{Bo1, Bo2}, and D. Bilyk, J. Dick, B. Green, L. Guth, and E. Saff for helpful conversations regarding point distributions on the sphere in Section \ref{sec:pt}. I am also grateful to C. Demeter, D. Jakobson, and S. Nonnenmacher for their encouragement during the preparation of this paper. 

My interest in spherical harmonics and their connection to quantum ergodicity began through discussions with S. Zelditch, whose insights greatly inspired the present research. He passed away in 2022, and I dedicate this paper to his memory.

I would like to thank the reviewer for the numerous suggestions, which have significantly improved the presentation of the article.

\section{Quantum ergodic spherical harmonics}\label{sec:QE}
In this section, we prove Theorem \ref{thm:QE}. Fix $\rho_0\in(0,1)$. Let $m\in\N$ satisfy $m\to\infty$ and $m=O(N^{\rho_0})$. Suppose that Gaussian beams $\{Q_j\}_{j=1}^m\subset\SH_N$ are chosen so that their poles $\{p_j\}_{j=1}^m$ satisfy the following conditions:
\begin{itemize}
\item[(S).] Separation \eqref{eq:S}: since $m=O(N^{\rho_0})$, after decreasing $c$ if necessary,
$$\dist\left(p_j,p_k\right)\ge\frac{c}{\sqrt m}\ge cN^{-\frac{\rho_0}{2}}\quad\text{for all }j\ne k,$$
\item[(E).] Equidistribution \eqref{eq:E}:
$$\lim_{m\to\infty}\frac1m\sum_{j=1}^m f\left(p_j\right)=\frac{1}{4\pi}\int_\Sb f(x)\,dx\quad\text{for any }f\in C^\infty(\Sb).$$
\end{itemize}

Throughout this section, we use both $N$ and the semiclassical parameter
$$h=h_N=\frac{1}{\sqrt{N(N+1)}}\approx N^{-1}$$
from Section \ref{sec:QE0}, so that $h\to0$ if and only if $N\to\infty$.

Define
$$u_N=\frac{F_N}{\left\|F_N\right\|_{L^2(\Sb)}},\qquad\text{where }F_N=\sum_{j=1}^mQ_j.$$
We show that $u_N$ are quantum ergodic; that is, for each $a\in C^\infty_0(T^*\Sb)$,
\begin{eqnarray}
&&\left\langle\Op_h(a)u_N,u_N\right\rangle\nonumber\\
&=&\frac{1}{\left\|F_N\right\|_{L^2(\Sb)}^2}\sum_{j,k=1,j\ne k}^m\left\langle\Op_h(a)Q_j,Q_k\right\rangle_{L^2(\Sb)}+\frac{1}{\left\|F_N\right\|_{L^2(\Sb)}^2}\sum_{j=1}^m\left\langle\Op_h(a)Q_j,Q_j\right\rangle_{L^2(\Sb)}\label{eq:OpuN}\\
&\to&\int_{S^*\Sb}a\,d\mu_L,\nonumber
\end{eqnarray}
where $\mu_L$ denotes the normalized Liouville measure on $S^*\Sb$.  

To prove this, in Section \ref{sec:FL} we estimate $\|F_N\|_{L^2(\Sb)}$. Section \ref{sec:SA} reviews the semiclassical analysis tools required for our argument. These are then applied in Section \ref{sec:matrix} to control the matrix elements $\langle\Op_h(a)Q_j,Q_k\rangle$. Finally, in Section \ref{sec:circle} we complete the proof by invoking a geometric representation of the cosphere bundle $S^*\Sb$ as the space of oriented great circles on $\Sb$.

\subsection{$L^2$ estimate}\label{sec:FL}
The $L^2$ bound for $F_N$ follows directly from the separation condition (S).  
By Han \cite[Lemma 5]{Ha1}, the inner product between two Gaussian beams $Q_j$ and $Q_k$ satisfies
\begin{equation}\label{eq:inner}
\left|\left\langle Q_j,Q_k\right\rangle\right|\le\left(\cos\frac{\dist\left(p_j,p_k\right)}{2}\right)^{2N}.
\end{equation}
Hence, the inner product decreases as the angular separation between $p_j$ and $p_k$ increases. In the extreme case $\dist(p_j,p_k)=\pi$, i.e., when the poles are antipodal, the corresponding beams are exactly orthogonal.

Since each $Q_j$ is $L^2$-normalized as in \eqref{eq:Q0},
\begin{eqnarray*}
\left\|F_N\right\|_{L^2(\Sb)}^2&=&\sum_{k,j=1}^m\left\langle Q_j,Q_k\right\rangle_{L^2(\Sb)}\\
&=&\sum_{j=1}^m\left\|Q_j\right\|_{L^2(\Sb)}^2+\sum_{k,j=1,\,j\ne k}^m\left\langle Q_j,Q_k\right\rangle_{L^2(\Sb)}\\
&=&m+\sum_{k,j=1,\,j\ne k}^m\left\langle Q_j,Q_k\right\rangle_{L^2(\Sb)}.
\end{eqnarray*}
Recall the standard estimate
\begin{equation}\label{eq:cos}
\cos\beta=1-\frac{\beta^2}{2}+O\left(\beta^4\right)\le1-\frac{\beta^2}{3}\quad\text{for }0\le\beta\le\frac\pi2.
\end{equation}
Since $\dist\left(p_j,p_k\right)\ge cN^{-\frac{\rho_0}{2}}$ for all $j\ne k$, we obtain from \eqref{eq:inner} and \eqref{eq:cos}
\begin{eqnarray*}
&&\sum_{k,j=1,\,j\ne k}^m\left|\left\langle Q_j,Q_k\right\rangle_{L^2(\Sb)}\right|\\
&\le&\sum_{k,j=1,\,j\ne k,\,\dist\left(p_j,p_k\right)\le\frac12}^m\left|\left\langle Q_j,Q_k\right\rangle_{L^2(\Sb)}\right|
+\sum_{k,j=1,\,\dist\left(p_j,p_k\right)>\frac12}^m\left|\left\langle Q_j,Q_k\right\rangle_{L^2(\Sb)}\right|\\
&\le&m^2\left[\cos\left(\frac c2N^{-\frac{\rho_0}{2}}\right)\right]^{2N}+m^2\left(\cos\frac14\right)^{2N}\\
&\le&O\left(N^{2\rho_0}\right)\cdot\exp\left[2N\log\left(1-\frac{c^2N^{-\rho_0}}{12}\right)\right]+O\left(N^{2\rho_0}\right)\cdot O\left(N^{-\infty}\right)\\
&\le&O\left(N^{2\rho_0}\right)\cdot\exp\left[-\frac{c^2N^{1-\rho_0}}{6}\right]+O\left(N^{-\infty}\right)\\
&=&O\left(N^{-\infty}\right)\quad\text{as }N\to\infty,
\end{eqnarray*}
because $0<\rho_0<1$. Hence,
\begin{equation}\label{eq:FL}
\left\|F_N\right\|_{L^2(\Sb)}^2=m+\sum_{k,j=1,\,j\ne k}^m\left\langle Q_j,Q_k\right\rangle_{L^2(\Sb)}=m+O\left(N^{-\infty}\right)\quad\text{as }N\to\infty.
\end{equation}

\subsection{Semiclassical preliminaries}\label{sec:SA}
We recall several standard facts from semiclassical analysis; see Zworski \cite{Zw} for a comprehensive exposition. Throughout, we consider only symbols with compact support, which suffices for our argument.

\begin{defn}[Symbol classes]
Let $\rho\in[0,\frac12)$. Fix $h_0\in(0,1)$. We say that $a(x,\xi;h)$ belongs to the symbol class $S_\rho^\comp(\M)$ if $a(x,\xi;h)\in C_0^\infty(T^*\M)$ for each $h\in(0,h_0)$ and for every multi-index pair $\alpha,\beta$,
$$\sup_{x\in\M,\xi\in T_x^*\M}\left|\partial^\alpha_x\partial^\beta_\xi a\right|\le Ch^{-\rho(|\alpha|+|\beta|)},$$
where $C=C(\alpha,\beta)>0$ is independent of $h$. The infimum $C(\alpha,\beta)$ for which these estimates hold is referred to as the seminorm of $a$. 
\end{defn}

If $a(x,\xi)\in C_0^\infty(T^*\M)$ is independent of $h$, then clearly $a\in S_0^\comp(\M)$.
  
We associate symbols in $S_\rho^\comp(\M)$ with semiclassical pseudodifferential operators as follows.  
For $a\in S_\rho^\comp(\R^d)$, define the left-quantization by
\begin{equation}\label{eq:leftq}
\Op_h(a)u(x)=\frac{1}{(2\pi h)^d}\int_{\R^d}\int_{\R^d}e^{\frac ih\langle\xi,x-y\rangle}a(x,\xi)u(y)\,d\xi dy,
\end{equation}
where $u\in C_0^\infty(\R^d)$. For $a\in S_\rho^\comp(\M)$, the operator $\Op_h(a)$ is defined locally using charts on $\M$.  
The correspondence between operators and symbols depends on the chosen quantization rule \eqref{eq:leftq} and is not unique. Nevertheless, the following standard results hold; see Zworski \cite[Sections 4.4 and 4.5]{Zw}.

$\bullet$ [$L^2$ boundedness]  
Let $\rho\in[0,\frac12)$ and $a\in S_\rho^\comp(\M)$. Then 
\begin{equation}\label{eq:OpL2}
\left\|\Op_h(a)\right\|_{L^2(\M)\to L^2(\M)}\le C,
\end{equation} 
where $C=C(a)>0$ depends on finitely many seminorms of $a$ and is independent of $h$.

$\bullet$ [Adjoint] Let $\rho\in[0,\frac12)$ and $a\in S_\rho^\comp(\M)$. Denote by $\Op_h(a)^\star$ the adjoint of $\Op_h(a)$ in $L^2(\M)$. Then
\begin{equation}\label{eq:adjoint}
\Op_h(a)^\star=\Op_h\left(\ol a\right)+O_{L^2(\M)\to L^2(\M)}\left(h^{1-2\rho}\right).
\end{equation}

$\bullet$ [Product] Let $\rho\in[0,\frac12)$ and $a,b\in S_\rho^\comp(\M)$. Then $ab\in S_\rho^\comp(\M)$ and
\begin{equation}\label{eq:prod}
\Op_h(a)\Op_h(b)=\Op_h(ab)+O_{L^2(\M)\to L^2(\M)}\left(h^{1-2\rho}\right).
\end{equation}

$\bullet$ [Microlocalization] Let $(h^2\Delta_\M+1)u_h=0$. Suppose that $a\in C^\infty_0(T^*\M)$ satisfies $a=1$ in a neighborhood of $S^*\M$. Then 
\begin{equation}\label{eq:micro}
\left\|\Op_h(a)u_h-u_h\right\|_{L^2(\M)}=O(h^\infty).
\end{equation}

The constants implied in \eqref{eq:adjoint}, \eqref{eq:prod}, and \eqref{eq:micro} depend on finitely many seminorms of the symbols involved and are independent of $h$.

\subsection{Estimates of the matrix elements}\label{sec:matrix}
We begin by estimating the diagonal terms in \eqref{eq:OpuN}, which correspond to the semiclassical measures associated with Gaussian beams.

\begin{prop}[Semiclassical measure of Gaussian beams]\label{prop:smG}
Let $Q_p$ be the Gaussian beam with pole $p\in\Sb$. Then
$$\left\langle\Op_h(a)Q_p,Q_p\right\rangle=\frac{1}{2\pi}\int_{G_p}a\left(x,\xi_p\right)\,dl+O(h)\quad\text{for all }a\in C_0^\infty(T^*\Sb),$$
where $\xi_p\in S^*G_p$ is the covector determined by the right-hand rule with respect to $p$, and $dl$ denotes the arclength measure on $G_p$ (see Figure \ref{fig:Gn}).
\end{prop}

\begin{figure}[H]
\begin{tikzpicture}
\draw[thick] (0,0) circle (4);

\draw[thick] (-4,0) arc (180:360:4 and 0.6);
\draw[thick,dashed] (4,0) arc (0:180:4 and 0.6);

\filldraw (0,4) circle (0.05cm) node[above]{The north pole $p_0$};
\draw (4,0) node[right]{The equator $G_0$};
\draw[very thick,->] (0,-0.6) -- (1,-0.6) node[midway,below]{$\xi_0$};
\draw[very thick,->] (0,0.6) -- (-1,0.6) node[midway,below]{$\xi_0$};
\end{tikzpicture}
\caption{The semiclassical measure of the Gaussian beam $Q_0$ is the normalized uniform measure supported on $G_0\times\{\xi_0\}\subset S^*\Sb$.}
\label{fig:Gn}
\end{figure}

\begin{proof}
This is a particular case of the semiclassical measure computed by Zelditch \cite[Proposition 12.1]{Ze3}; see also Zworski \cite[Example 2 in Section 5.1]{Zw}. The normalization factor $\frac{1}{2\pi}$ can be verified directly by setting $a=1$ in a neighborhood of $G_p$.
\end{proof}

We next analyze the off-diagonal terms $\langle\Op_h(a)Q_j,Q_k\rangle$ for $j\ne k$ in \eqref{eq:OpuN}, outlining the main idea. By Proposition \ref{prop:smG}, each $Q_j$ is microlocalized near $\xi_j$ in frequency space, while $Q_k$ is microlocalized near $\xi_k$. By the separation condition (S), the angle between $p_j$ and $p_k$ satisfies
$$\beta=\dist\left(p_j,p_k\right)\ge cN^{-\frac{\rho_0}{2}}\ge ch^{\frac{\rho_0}{2}}.$$
The same separation holds for $\xi_j$ and $\xi_k$, implying that $Q_j$ and $Q_k$ are microlocalized in  disjoint regions separated by $ch^{\frac{\rho_0}{2}}\gg h^\frac12$, because $\rho_0<1$. Since a semiclassical pseudodifferential operator $\Op_h(a)$ preserves microlocal support, it follows that $\langle\Op_h(a)Q_j,Q_k\rangle$ is negligible for $j\ne k$. We now state this estimate, which is the conclusion of the analysis of the present section.

\begin{lemma}[Off-diagonal estimate]\label{lemma:off}
Let $0<\rho_0<1$ and $a\in C_0^\infty(T^*\Sb)$. Let $Q_j,Q_k\in\SH_N$ be Gaussian beams whose poles satisfy
$$\dist\left(p_j,p_k\right)\ge cN^{-\frac{\rho_0}{2}}$$
for some $c>0$. Then
$$\left\langle\Op_h(a)Q_j,Q_k\right\rangle_{L^2(\Sb)}=O\left(h^\infty\right)\quad\text{as }h\to0,$$
where the implicit constants depend only on $\rho_0$, $c$, and finitely many seminorms of $a$.
\end{lemma}

\begin{proof}
Without loss of generality, we take $Q_j=Q_0$. In spherical coordinates $(\phi,\theta)$, let $\chi_0\in C_0^\infty(\Sb)$ be independent of $\theta$, with $0\le\chi_0\le1$ and
$$\chi_0(\phi)=\begin{cases}
1, & \text{if }\left|\phi-\frac\pi2\right|\le\frac12h^{\frac{\rho_0}{2}},\\ 
0, & \text{if }\left|\phi-\frac\pi2\right|\ge h^{\frac{\rho_0}{2}}.
\end{cases}$$
If $(\phi,\theta)\notin\supp\chi_0$, then 
$$|\alpha|=\left|\frac\pi2-\phi\right|\ge h^{\frac{\rho_0}{2}}.$$
Using \eqref{eq:Q0polar} and \eqref{eq:cos},
$$\left|Q_0(\phi,\theta)\right|\le c_2N^\frac14(\cos\alpha)^N\le c_2N^\frac14\exp\left[N\log\left(1-\frac{h^{\rho_0}}{3}\right)\right]\le cN^\frac14\exp\left[-\frac{h^{\rho_0-1}}{3}\right]=O\left(h^\infty\right),$$
since $\rho_0-1<0$ and $h=\frac{1}{\sqrt{N(N+1)}}\approx N^{-1}$. Hence,
$$Q_0=\chi_0Q_0+O\left(h^\infty\right)\quad\text{uniformly on }\Sb.$$
Similarly, define $\chi_k\in C_0^\infty(\Sb)$ supported in the $h^{\frac{\rho_0}{2}}$-neighborhood of the great circle $G_k$. Then
$$Q_k=\chi_kQ_k+O\left(h^\infty\right)\quad\text{uniformly on }\Sb.$$
Using the $L^2$-boundedness of $\Op_h(a)$ from \eqref{eq:OpL2},
$$\left\langle\Op_h(a)Q_0,Q_k\right\rangle_{L^2(\Sb)}=\left\langle\Op_h(a)\left(\chi_0Q_0\right),\chi_kQ_k\right\rangle_{L^2(\Sb)}+O\left(h^\infty\right).$$
Therefore, it suffices to estimate
$$\Op_h(a)\left(\chi_0(\tilde\phi)Q_0(\tilde\phi,\tilde\theta)\right)(\phi,\theta)\quad\text{and}\quad Q_k(\phi,\theta)$$
for $(\tilde\phi,\tilde\theta)\in\Nc_{h^{\frac{\rho_0}{2}}}(G_0)$ and $(\phi,\theta)\in\Nc_{h^{\frac{\rho_0}{2}}}(G_k)$, see Figure \ref{fig:inter}.  

\begin{figure}[H]
\begin{tikzpicture}
\draw[thick] (0,0) circle (4);

\filldraw (0,4) circle (0.05cm) node[above]{The north pole $p_0$};
\draw (4,0) node[right]{The equator $G_0$};
\draw[thick] (0,4) -- (0,0);
\draw[thick] (-4,0) arc (180:360:4 and 0.6);
\draw[thick,dashed] (-4,0) arc (180:0:4 and 0.6);
\draw[] (-3.967,0.5) arc (180:360:3.967 and 0.6);
\draw[dashed] (3.967,0.5) arc (0:180:3.967 and 0.6);
\draw[] (-3.967,-0.5) arc (180:360:3.967 and 0.6);
\draw[dashed] (3.967,-0.5) arc (0:180:3.967 and 0.6);
\draw[thick,<->] (-4.2,-0.5) arc (190:175:4.2) node[left,midway]{$2h^{\frac{\rho_0}{2}}$};

\filldraw[rotate around={0.54r:(0,0)}] (0,4) circle (0.05cm) node[anchor=south east]{The pole $p_k$ of $G_k$};
\draw (3.708,2) node[right]{A great circle $G_k$};
\draw[rotate around={0.54r:(0,0)},thick] (0,4) -- (0,0);
\draw[rotate around={0.54r:(0,0)},thick] (-4,0) arc (180:360:4 and 0.6);
\draw[rotate around={0.54r:(0,0)},thick,dashed] (-4,0) arc (180:0:4 and 0.6);
\draw[rotate around={0.54r:(0,0)}] (-3.967,0.5) arc (180:360:3.967 and 0.6);
\draw[rotate around={0.54r:(0,0)},dashed] (3.967,0.5) arc (0:180:3.967 and 0.6);
\draw[rotate around={0.54r:(0,0)}] (-3.967,-0.5) arc (180:360:3.967 and 0.6);
\draw[rotate around={0.54r:(0,0)},dashed] (3.967,-0.5) arc (0:180:3.967 and 0.6);
\draw[thick,<->] (-3.9,-1.7) arc (203:219:4) node[anchor=north east,midway]{$2h^{\frac{\rho_0}{2}}$};

\draw[thick] (0,2.2) arc (90:120:2.2) node[midway,above]{$\beta$};
\end{tikzpicture}
\caption{Intersection of $\Nc_{h^{\frac{\rho_0}{2}}}(G_k)$ and $\Nc_{h^{\frac{\rho_0}{2}}}(G_0)$.}
\label{fig:inter}
\end{figure}

In $\Nc_{h^{\frac{\rho_0}{2}}}(G_k)$,
$$Q_k(\phi,\theta)=C_NN^\frac14\left(\sin(\phi+\beta)\right)^Ne^{iN\theta\cos\beta}.$$
In local coordinates,
\begin{eqnarray*}
&&\Op_h(a)\left(\chi_0Q_0\right)(\phi,\theta)\\
&=&\frac{1}{(2\pi h)^2}\int_{\R^2}\int_{\R^2}e^{\frac ih(\langle\xi_\phi,\phi-\tilde\phi\rangle+\langle\xi_\theta,\theta-\tilde\theta\rangle)}a\left(\phi,\theta,\xi_\phi,\xi_\theta\right)\chi_0(\tilde\phi)Q_0(\tilde\phi,\tilde\theta)\,d\xi_\phi d\xi_\theta d\tilde\phi d\tilde\theta\\
&=&\frac{C_NN^\frac14}{(2\pi h)^2}\int_{\R^2}\int_{\R^2}e^{\frac ih(\langle\xi_\phi,\phi-\tilde\phi\rangle+\langle\xi_\theta,\theta-\tilde\theta\rangle)}a\left(\phi,\theta,\xi_\phi,\xi_\theta\right)\chi_0(\tilde\phi)\left(\sin\tilde\phi\right)^Ne^{iN\tilde\theta}\,d\xi_\phi d\xi_\theta d\tilde\phi d\tilde\theta.
\end{eqnarray*}
Hence,
\begin{eqnarray*}
&&\left\langle\Op_h(a)\left(\chi_0Q_0\right),\chi_kQ_k\right\rangle_{L^2(\Sb)}\\
&=&\frac{C_N^2N^\frac12}{(2\pi h)^2}\int_{\R^2}\int_{\R^2}\int_{\R^2}e^{\frac ih(\langle\xi_\phi,\phi-\tilde\phi\rangle+\langle\xi_\theta,\theta-\tilde\theta\rangle)}a\left(\phi,\theta,\xi_\phi,\xi_\theta\right)\chi_0(\tilde\phi)\left(\sin\tilde\phi\right)^Ne^{iN\tilde\theta}\\
&&\cdot\chi_k(\phi,\theta)\left(\sin(\phi+\beta)\right)^Ne^{-iN\theta\cos\beta}\,d\xi_\phi d\xi_\theta d\tilde\phi d\tilde\theta d\phi d\theta\\
&=&\frac{C_N^2h^{-\frac12}}{(2\pi h)^2}\int_{\R^2}\int_{\R^2}e^{\frac ih\langle\xi_\phi,\phi-\tilde\phi\rangle}\chi_0(\tilde\phi)(\sin\phi)^\frac1h\left(\sin(\phi+\beta)\right)^\frac1h\\
&&\left[\int_\R a\left(\phi,\theta,\xi_\phi,\xi_\theta\right)\chi_k(\phi,\theta)e^{\frac ih(1-\cos\beta)\theta}\,d\theta\right]\,d\xi_\phi d\xi_\theta d\tilde\phi  d\phi+O\left(h^\infty\right).
\end{eqnarray*}
The phase function of the $\theta$-integral is
$$\Phi(\theta)=(1-\cos\beta)\theta.$$
Since $\beta\ge ch^{\frac{\rho_0}{2}}$, by \eqref{eq:cos}, 
$$|\Phi'(\theta)|=|1-\cos\beta|\ge\frac{\beta^2}{3}\ge ch^{\rho_0}.$$
As $\rho_0<1$, repeated integration by parts in $\theta$ applies, giving
$$\int_\R a\left(\phi,\theta,\xi_\phi,\xi_\theta\right)\chi_k(\phi,\theta)e^{\frac ih(1-\cos\beta)\theta}\,d\theta=O\left(h^\infty\right),$$
see Zworski \cite[Lemma 3.14]{Zw}. Consequently,
$$\left\langle\Op_h(a)Q_0,Q_k\right\rangle_{L^2(\Sb)}=O\left(h^\infty\right).$$
The same estimate holds for $\langle\Op_h(a)Q_j,Q_k\rangle_{L^2(\Sb)}$ whenever $j\ne k$: applying a rotation $R\in\SO(3)$ with $Rp_j=p_0$ replaces $a$ by its pullback under $R$, whose seminorms are comparable to those of $a$, and leaves the separation hypothesis unchanged.
\end{proof}

Combining this with the $L^2$ estimate of $F_N$ in \eqref{eq:FL} and the diagonal estimate in Proposition \ref{prop:smG}, we obtain from \eqref{eq:OpuN}:
\begin{eqnarray*}
&&\left\langle\Op_h(a)u_N,u_N\right\rangle\\
&=&\frac{1}{\left\|F_N\right\|_{L^2(\Sb)}^2}\sum_{j,k=1,j\ne k}^m\left\langle\Op_h(a)Q_j,Q_k\right\rangle_{L^2(\Sb)}+\frac{1}{\left\|F_N\right\|_{L^2(\Sb)}^2}\sum_{j=1}^m\left\langle\Op_h(a)Q_j,Q_j\right\rangle_{L^2(\Sb)}\\
&=&\frac{1}{m+O\left(h^\infty\right)}\cdot m^2\cdot O\left(h^\infty\right)+\frac{1}{m+O\left(h^\infty\right)}\sum_{j=1}^m\left(\frac{1}{2\pi}\int_{G_j}a\left(x,\xi_j\right)\,dl+O(h)\right)\\
&=&\frac{1}{2\pi m}\sum_{j=1}^m\int_{G_j}a\left(x,\xi_j\right)\,dl+O(h),
\end{eqnarray*}
where $\xi_j\in S^*G_j$ is determined by the right-hand rule with respect to $p_j$.

\subsection{The space of oriented great circles}\label{sec:circle}
The equidistribution condition (E) \eqref{eq:E} concerns smooth functions of points on $\Sb$, whereas the quantities appearing at the end of Section \ref{sec:matrix} are averages of $a$ over great circles in $S^*\Sb$, parametrized by their poles. The following lemma reconciles the two.

\begin{lemma}\label{lemma:avg}
For $a\in C_0^\infty(T^*\Sb)$,
$$A_a(p)=\frac{1}{2\pi}\int_{G_p}a\left(x,\xi_p\right)\,dl,\quad\text{where }p\in\Sb,$$
defines a function $A_a\in C^\infty(\Sb)$. 
\end{lemma}

\begin{proof}
The equator $G_0$, with the north pole $p_0$ as its pole, has the arclength parametrization $\gamma_0(t)=(\cos t,\sin t,0)$ with $t\in[0,2\pi)$, so that $\gamma_0'(t)=\xi_0$ at $\gamma_0(t)$ and
$$A_a\left(p_0\right)=\frac{1}{2\pi}\int_0^{2\pi}a\left(\gamma_0(t),\gamma_0'(t)\right)\,dt,$$
where we identify the tangent vector and the cotangent vector on $\Sb$ by the round metric.

Let $R\in\SO(3)$ be a rotation with $Rp_0=p$. Since $R$ is an orientation-preserving isometry, it maps $\gamma_0$ to the arclength parametrization $\gamma_p=R\gamma_0$ of $G_p$ and carries $\xi_0$ to $\xi_p$. It then follows that
$$A_a(p)=\frac{1}{2\pi}\int_0^{2\pi}a\left(\gamma_p(t),\gamma'_p(t)\right)\,dt.$$
Two rotations taking $p_0$ to $p$ differ by a rotation about $p_0$, which amounts to a translation in $t$; as the integral is taken over a full period, the right-hand side is independent of the choice of $R$. Choosing $R=R(p)$ to depend smoothly on $p$, which is possible near any point of $\Sb$ because $R\mapsto Rp_0$ is a submersion from $\SO(3)$ onto $\Sb$, and differentiating under the integral sign, we conclude that $A_a\in C^\infty(\Sb)$.
\end{proof}

Since $G_j=G_{p_j}$ and $\xi_j=\xi_{p_j}$, the estimate at the end of Section \ref{sec:matrix} reads
$$\left\langle\Op_h(a)u_N,u_N\right\rangle=\frac1m\sum_{j=1}^mA_a\left(p_j\right)+O(h).$$
Applying the equidistribution condition (E) \eqref{eq:E} to the function $f=A_a\in C^\infty(\Sb)$ of Lemma \ref{lemma:avg}, and recalling that $m=m(N)\to\infty$ as $N\to\infty$, we obtain
\begin{eqnarray*}
\left\langle\Op_h(a)u_N,u_N\right\rangle&=&\frac1m\sum_{j=1}^mA_a\left(p_j\right)+O(h)\\
&\to&\frac{1}{4\pi}\int_\Sb A_a(p)\,dp\\
&=&\frac{1}{8\pi^2}\int_\Sb\int_{G_p}a\left(x,\xi_p\right)\,dldp.
\end{eqnarray*}
To conclude the proof of Theorem \ref{thm:QE}, recall that the cosphere bundle $S^*\Sb$ can be identified with the space of oriented great circles, each parametrized by its pole; see Jakobson–Zelditch \cite[Section 3]{JZ}. Hence,
$$\frac{1}{8\pi^2}\int_{\Sb}\int_{G_p}a\left(x,\xi_p\right)\,dldp=\int_{S^*\Sb}a\,d\mu_L.$$
The constant $\frac{1}{8\pi^2}$ corresponds to the normalization $\mu_L(S^*\Sb)=1$. Indeed, taking $a\in C^\infty_0(T^*\Sb)$ with $a=1$ in a neighborhood of $S^*\Sb$, we have by \eqref{eq:micro},
$$\langle\Op_h(a)u_N,u_N\rangle_{L^2(\Sb)}=1+O(h^\infty),$$
confirming that both sides agree and establishing the normalization.

\section{Spherical harmonics with controlled $L^\infty$ norms}\label{sec:BddSHs}
In this section, we prove Theorem \ref{thm:BddSHs}. Let $0<\delta<1$ be a small parameter to be specified later. Set
\begin{equation}\label{eq:m}
m=\left\lfloor\delta N^{\frac12}\right\rfloor+1.
\end{equation}
Suppose that Gaussian beams $\{Q_j\}_{j=1}^m\subset\SH_N$ are chosen so that their poles $\{p_j\}_{j=1}^m$ satisfy the following conditions:
\begin{itemize}
\item[(S).] Separation \eqref{eq:S}: since $m\le2N^\frac12$, after decreasing $c$ if necessary,
$$\dist\left(p_j,p_k\right)\ge\frac{c}{\sqrt m}\ge cN^{-\frac14}\quad\text{for all }j\ne k.$$
\item[(G).] Bounded clustering around great circles \eqref{eq:G}: there exists $L(m)\ge1$ such that
$$\#\left\{j:p_j\in\Nc_{m^{-1}}(G)\right\}\le L(m)\quad\text{for all great circles }G\subset\Sb.$$
\end{itemize}
Define
$$u_N=\frac{F_N}{\left\|F_N\right\|_{L^2(\Sb)}},\qquad\text{where }F_N=\sum_{j=1}^mQ_j.$$
We shall show that
$$\left\|u_N\right\|_{L^\infty(\Sb)}\le C\delta^{-\frac12}L(m),$$
for some absolute constant $C>0$.

To this end, in Section \ref{sec:FL2} we establish
\begin{equation}\label{eq:FL2}
\left\|F_N\right\|_{L^2(\Sb)}=\sqrt m+O\left(N^{-\infty}\right),
\end{equation}
and in Section \ref{sec:FLinf} we prove
\begin{equation}\label{eq:FLinf}
\left\|F_N\right\|_{L^\infty(\Sb)}\le cN^\frac14L(m),
\end{equation}
where $c>0$ is an absolute constant. Combining these two estimates gives
$$\left\|u_N\right\|_{L^\infty(\Sb)}\le\frac{cN^\frac14L(m)}{\sqrt m+O\left(N^{-\infty}\right)}\le C\delta^{-\frac12}L(m).$$
Hence, Theorem \ref{thm:BddSHs} follows.

\subsection{$L^2$ estimate}\label{sec:FL2}
Since $m\le2N^\frac12=O(N^\frac12)$, the $L^2$ bound for $F_N$ follows directly from the separation condition (S), in the same way as in Section \ref{sec:FL}. We therefore omit the proof of \eqref{eq:FL2}.

\subsection{$L^\infty$ estimate}\label{sec:FLinf}
The $L^\infty$ bound for $F_N$ follows from the bounded clustering condition (G). Without loss of generality, we estimate $F_N(p_0)$ at the north pole $p_0$ (since one may always rotate $F_N$ so that its maximum is attained at $p_0$). By the coordinate expression \eqref{eq:Q0polar} for Gaussian beams,
$$\left|F_N\left(p_0\right)\right|=\left|\sum_{j=1}^mQ_j\left(p_0\right)\right|\le\sum_{j=1}^m\left|Q_j\left(p_0\right)\right|=C_NN^\frac14\sum_{j=1}^m\left|\sin\phi_j\right|^N\le c_2N^\frac14\sum_{j=1}^m\left(\cos\alpha_j\right)^N,$$
where $\phi_j=\dist(p_j,p_0)$ and $\alpha_j=\frac\pi2-\phi_j$, i.e., the angle between $p_j$ and the equator $G_0$. 

Partition $\Sb$, using our coordinates, by the latitudes
$$\left\{\alpha=lm^{-1}\right\},\quad l=0,\pm1,\dots,\pm\left(\left\lfloor\frac{\pi m}{2}\right\rfloor+1\right),$$
and denote the $l$-th strip by
$$S_l=\left\{lm^{-1}\le\alpha<(l+1)m^{-1}\right\}.$$
The poles $\{p_j\}_{j=1}^m$ then fall into three groups:

\begin{enumerate}[(I).]
\item $\Gc_1=\{j:|\alpha_j|\le m^{-1}\}$: poles within the $m^{-1}$-neighborhood of the equator $G_0$.  
By Condition (G), the number of such poles is bounded by $L(m)$.  

\item $\Gc_2=\{j:m^{-1}<|\alpha_j|\le\frac13\}$: poles lying in intermediate latitudinal strips.  
Each strip $S_l$ can be covered by $O(l)$ neighborhoods $\Nc_{m^{-1}}(G)$ of great circles, hence contains at most $O(l)\cdot L(m)$ poles; see Lemma \ref{lemma:cover}.  

\item $\Gc_3=\{j:|\alpha_j|>\frac13\}$: poles located near the north and south poles of $\Sb$, whose contributions decay exponentially in $N$.  
\end{enumerate}

We now estimate the contributions from each group.

\textbf{Group I.} Applying the bounded clustering condition (G) to $G_0$ gives
$$\#\left\{j:\left|\alpha_j\right|\le\frac1m\right\}=\#\left\{j:p_j\in\Nc_{m^{-1}}\left(G_0\right)\right\}\le L(m),$$
and therefore
$$\sum_{j\in\Gc_1}\left(\cos\alpha_j\right)^N\le L(m).$$

\textbf{Group II.} For the $l$-th strip $S_l$ with $l\ge1$ (the case $l\le-1$ is symmetric), we first show that $S_l$ is covered by $O(l)$ neighborhoods of great circles.

\begin{lemma}[Covering of latitudinal strips]\label{lemma:cover}
Let $m\ge6$ and $1\le l\le\frac m3$. Then $S_l$ is contained in the union of the $m^{-1}$-neighborhoods $\Nc_{m^{-1}}(G)$ of at most $5(l+1)$ great circles $G\subset\Sb$.
\end{lemma}

\begin{proof}
We identify the points of $\Sb$ with unit vectors in $\R^3$, so that $\langle x,y\rangle=\cos\dist(x,y)$ for $x,y\in\Sb$. Since the great circle $G_p$ with pole $p$ consists of the points at distance $\frac\pi2$ from $p$,
$$\dist\left(x,G_p\right)=\left|\frac\pi2-\dist(x,p)\right|=\left|\arcsin\left\langle x,p\right\rangle\right|\quad\text{for }x\in\Sb.$$
Because $m^{-1}\le\frac\pi2$, it follows that
\begin{equation}\label{eq:nbhd}
x\in\Nc_{m^{-1}}\left(G_p\right)\quad\text{if and only if}\quad\left|\left\langle x,p\right\rangle\right|\le\sin\frac1m.
\end{equation}
Let $G$ be the great circle forming the angle $(l+1)m^{-1}$ with the equator $G_0$ whose highest point has longitude $0$; see Figure \ref{fig:strip}. Its pole is then
$$p=\left(-\sin\frac{l+1}{m},0,\cos\frac{l+1}{m}\right),$$
and
$$\frac{l+1}{m}\le\frac13+\frac16=\frac12,$$
since $l\le\frac m3$ and $m\ge6$.

Let $x\in S_l$ have latitude $\alpha$ and longitude $\theta$ with $|\theta|\le(l+1)^{-1}$. In our coordinates,
\begin{eqnarray*}
\left\langle x,p\right\rangle&=&-\cos\alpha\cos\theta\sin\left(\frac{l+1}{m}\right)+\sin\alpha\cos\left(\frac{l+1}{m}\right)\\
&=&-\sin\left(\frac{l+1}{m}-\alpha\right)+\cos\alpha\sin\left(\frac{l+1}{m}\right)\left(1-\cos\theta\right).
\end{eqnarray*}
Since $lm^{-1}\le\alpha<(l+1)m^{-1}$ for $x\in S_l$, the first term satisfies
$$0<\sin\left(\frac{l+1}{m}-\alpha\right)\le\sin\frac1m.$$
Moreover, by $\sin\beta\le\beta$ and $1-\cos\beta\le\frac{\beta^2}{2}$,
$$0\le\cos\alpha\sin\frac{l+1}{m}\left(1-\cos\theta\right)\le\frac{l+1}{m}\cdot\frac{1}{2(l+1)^2}\le\frac{1}{2m}\le\sin\frac1m,$$
where in the last inequality we used $\sin\beta\ge\frac{2\beta}{\pi}\ge\frac\beta2$ for $0\le\beta\le\frac\pi2$. Combining the last three displays,
$$\left|\left\langle x,p\right\rangle\right|\le\sin\frac1m,$$
and therefore $x\in\Nc_{m^{-1}}(G)$ by \eqref{eq:nbhd}.

Thus $\Nc_{m^{-1}}(G)$ contains all points of $S_l$ whose longitude lies within $(l+1)^{-1}$ of $0$. Set $K=\lceil\pi(l+1)\rceil$ and rotate $G$ about the polar axis by the angles $\frac{2\pi k}{K}$, $k=1,...,K$. Every longitude lies within $\frac\pi K\le\frac{1}{l+1}$ of $\frac{2\pi k}{K}$ modulo $2\pi$ for some $k$, while $S_l$ is invariant under rotations about the polar axis. Hence these great circles cover $S_l$, and their number is
$$K\le\pi(l+1)+1\le5(l+1).$$
\end{proof}

By Lemma \ref{lemma:cover} and Condition (G) applied to each of these great circles, the number of poles in $S_l$ is at most
$$5(l+1)L(m)\le10lL(m)\quad\text{for }1\le l\le\frac m3.$$

\begin{figure}[H]
\begin{tikzpicture}
\draw[thick] (0,0) circle (4);
\fill[] (0,0) circle (0.05cm) node[right]{$o$};
\draw[thick] (-3.708,1.5) arc (180:360:3.708 and 0.6);
\draw[thick,dashed] (3.708,1.5) arc (0:180:3.708 and 0.6);
\draw[thick] (-3.464,2) arc (180:360:3.464 and 0.6);
\draw[thick,dashed] (3.464,2) arc (0:180:3.464 and 0.6);
\draw[rotate around={0.54r:(0,0)},thick] (-4,0) arc (180:360:4 and 0.6) node[below,midway]{$G$};
\draw[rotate around={0.54r:(0,0)},thick,dashed] (-4,0) arc (180:0:4 and 0.6);
\draw[rotate around={0.54r:(0,0)}] (-3.967,0.5) arc (180:360:3.967 and 0.6);
\draw[rotate around={0.54r:(0,0)},dashed] (3.967,0.5) arc (0:180:3.967 and 0.6);
\draw[rotate around={0.54r:(0,0)}] (-3.967,-0.5) arc (180:360:3.967 and 0.6);
\draw[rotate around={0.54r:(0,0)},dashed] (3.967,-0.5) arc (0:180:3.967 and 0.6);
\draw[thick,<->] (-3.708,-1.5) arc (202:219:4) node[anchor=north east,midway]{$2m^{-1}$};
\draw (3.708,1.5) node[anchor=south west]{$S_l$};
\draw[] (-3.708,1.5) node[anchor=north east]{$lm^{-1}$};
\draw[] (-3.464,2) node[anchor=south east]{$(l+1)m^{-1}$};
\end{tikzpicture}
\caption{Intersection of a latitudinal strip $S_l$ and $\Nc_{m^{-1}}(G)$, where $G$ forms an angle of $(l+1)m^{-1}$ with the equator $G_0$.}
\label{fig:strip}
\end{figure}

If $p_j\in S_l$ with $1\le l\le\frac m3$, then $lm^{-1}\le\alpha_j\le\frac13<\frac\pi2$, so that $\cos\alpha_j\le\cos(lm^{-1})$. Using $m=\lfloor\delta N^\frac12\rfloor+1\le2\delta N^\frac12$ for $N$ sufficiently large, and \eqref{eq:cos},
\begin{eqnarray*}
\sum_{j\in\Gc_2}\left(\cos\alpha_j\right)^N
&\le&2\sum_{l=1}^{\left\lfloor\frac m3\right\rfloor}\sum_{j:\,p_j\in S_l}\left(\cos\alpha_j\right)^N\\
&\le&20L(m)\sum_{l=1}^\infty l\left[\cos\left(lm^{-1}\right)\right]^N\\
&\le&20L(m)\sum_{l=1}^\infty l\exp\left[N\log\left(1-\frac{l^2}{3m^2}\right)\right]\\
&\le&20L(m)\sum_{l=1}^\infty l\exp\left[-N\cdot\frac{l^2}{3\left(2\delta N^\frac12\right)^2}\right]\\
&\le&20L(m)\sum_{l=1}^\infty l\exp\left[-\frac{l^2}{12\delta^2}\right]\\
&\le&L(m)\quad\text{for sufficiently small }\delta,
\end{eqnarray*}
since the series in the last line tends to $0$ as $\delta\to0$.

\textbf{Group III.} Since the total number of poles is at most $m$,
$$\sum_{j\in\Gc_3}\left(\cos\alpha_j\right)^N\le m\left(\cos\frac13\right)^N\le2\delta N^\frac12\left(\cos\frac13\right)^N=O\left(N^{-\infty}\right)$$
as $N\to\infty$.

Combining the three contributions gives
$$\left|F_N\left(p_0\right)\right|\le c_2N^\frac14\sum_{j=1}^m\left(\cos\alpha_j\right)^N\le c_2N^\frac14\left[2L(m)+O\left(N^{-\infty}\right)\right]\le 3c_2N^\frac14L(m),$$
which proves \eqref{eq:FLinf}.

\appendix
\section{A point configuration by a probabilistic approach}\label{sec:prob}
We construct point configurations satisfying (S) separation \eqref{eq:S}, (E) equidistribution \eqref{eq:E}, and (G) bounded clustering around great circles \eqref{eq:G} with logarithmic $L(m)$. The argument is inspired by Demeter-Zhang \cite[Proposition 2.1]{DZ}, which establishes an analogous statement for point distributions on the plane; see also Carbery \cite{C} for an earlier result.

\begin{prop}\label{prop:prob}
There exist absolute constants $c,C>0$ and, for each sufficiently large $m$, a point configuration $\{p_j\}_{j=1}^m\subset\Sb$, such that the resulting family satisfies the following properties.
\begin{enumerate}
\item[(S).] [Separation] For all $j\ne k$,
$$\dist\left(p_j,p_k\right)\ge\frac{c}{\sqrt m}.$$
\item[(E).] [Equidistribution] For any $f\in C^\infty(\Sb)$,
$$\lim_{m\to\infty}\frac1m\sum_{j=1}^mf\left(p_j\right)=\frac{1}{4\pi}\int_\Sb f(x)\,dx.$$
\item[(G).] [Bounded clustering around great circles] For all great circles $G\subset\Sb$,
$$\#\left\{j:p_j\in\Nc_{m^{-1}}(G)\right\}\le\frac{C\log m}{\log\log m}.$$
\end{enumerate}
\end{prop}

\begin{proof}
For each $m$, fix a configuration $\{q_j\}_{j=1}^m$ so that the resulting family satisfies (E), and (S) with the constant $2c$; see Section \ref{sec:pt} for explicit constructions. Set
$$r=\frac{c}{\sqrt m}.$$
By (S), the geodesic balls $\{B(q_j,r)\}_{j=1}^m$ are pairwise disjoint. For each $j=1,...,m$, choose $p_j\in B(q_j,\frac r2)$ independently at random according to the normalized Riemannian area measure. 

We first verify (S) for $\{p_j\}_{j=1}^m$. For $j\ne k$, we have $\dist(q_j,q_k)\ge2cm^{-\frac12}=2r$ by (S) for $\{q_j\}_{j=1}^m$, while $\dist(p_j,q_j)\le\frac r2$ and $\dist(p_k,q_k)\le\frac r2$. Hence
$$\dist\left(p_j,p_k\right)\ge2r-\frac r2-\frac r2=r=\frac{c}{\sqrt m}.$$

Next, we verify (E) for $\{p_j\}_{j=1}^m$. For any $f\in C^\infty(\Sb)$,
\begin{eqnarray*}
\left|\frac1m\sum_{j=1}^mf\left(p_j\right)-\frac1m\sum_{j=1}^mf\left(q_j\right)\right|&\le&\frac1m\sum_{j=1}^m\left|f\left(p_j\right)-f\left(q_j\right)\right|\\
&\le&\frac1m\sum_{j=1}^m\|\nabla f\|_{L^\infty(\Sb)}\cdot\dist\left(p_j,q_j\right)\\
&\le&\|\nabla f\|_{L^\infty(\Sb)}\cdot\frac r2\\
&=&O_f\left(m^{-\frac12}\right)\to0\quad\text{as }m\to\infty.
\end{eqnarray*}
Thus (E) follows from the corresponding property of $\{q_j\}_{j=1}^m$.

We now show that (G) holds with high probability as $m\to\infty$. The probability space is
$$\prod_{j=1}^mB\left(q_j,\frac r2\right)\quad\text{equipped with the measure }\prod_{j=1}^m\frac{d\Area}{\Area\left(B\left(q_j,\frac r2\right)\right)}.$$
Fix a great circle $G\subset\Sb$. Consider its $2m^{-1}$-neighborhood,
$$\Nc_{2m^{-1}}(G)=\left\{x\in\Sb:\dist(G,x)\le\frac2m\right\}.$$ 
For $j=1,...,m$, define independent Bernoulli variables
$$X_j=\indicator_{\Nc_{2m^{-1}}(G)}\left(p_j\right).$$ 
Then 
$$X=\sum_{j=1}^m X_j$$
counts the number of points $p_j$ lying in $\Nc_{2m^{-1}}(G)$.

Since $\Area(B(q_j,\frac r2)\cap\Nc_{2m^{-1}}(G))\le brm^{-1}$, we have
$$\Pr\left(X_j=1\right)=\Pr\left(p_j\in\Nc_{2m^{-1}}(G)\right)\le\frac{\Area\left(B\left(q_j,\frac r2\right)\cap\Nc_{2m^{-1}}(G)\right)}{\Area\left(B\left(q_j,\frac r2\right)\right)}\le\frac{b}{\sqrt m}.$$
Here and throughout, $b$ denotes a positive absolute constant, which may take different values from line to line.

Because of (S) separation, $\Nc_{2m^{-1}}(G)$ can intersect at most $bm^\frac12$ balls among $\{B(q_j,r)\}_{j=1}^m$. Hence,
$$\E[X]=\sum_{j=1}^m\Pr\left(X_j=1\right)\le bm^\frac12\cdot\frac{b}{\sqrt m}\le b.$$
Therefore, the expected number of points in the $m^{-1}$-neighborhood of any fixed great circle is uniformly bounded.

\begin{figure}[H]
\begin{tikzpicture}
\draw[thick] (-3,0.5) circle (2);
\filldraw[thick,draw=black,fill=gray!20] (-3,0.5) circle (1);
\draw[->] (-3,0.5) -- (-3,2.5);
\draw[] (-3,2.25) node[right]{$r$};
\filldraw[black] (-3,0.5) circle (0.05) node[right]{$q_j$};
\filldraw[black] (-3,-0.15) circle (0.05) node[right]{$p_j$};

\draw[thick] (3,-0.5) circle (2);
\filldraw[thick,draw=black,fill=gray!20] (3,-0.5) circle (1);
\draw[->] (3,-0.5) -- (3,-2.5);
\draw[] (3,-2.25) node[right]{$r$};
\filldraw[black] (3,-0.5) circle (0.05) node[right]{$q_k$};
\filldraw[black] (3,-0.15) circle (0.05) node[right]{$p_k$};

\draw[thick,<->] (-6,-0.25) -- (-6,0.25) node[left,midway]{$4m^{-1}$};
\draw[thick] (-6,-0.25) -- (6,-0.25);
\draw[thick] (-6,0.25) -- (6,0.25);
\draw[thick,-] (6,-0.25) -- (6,0.25) node[right,midway]{$\Nc_{2m^{-1}}(G)$};
\end{tikzpicture}
\caption{Intersections of $\Nc_{2m^{-1}}(G)$ with $B(q_j,r)$}
\label{fig:int}
\end{figure}

We apply a standard large-deviation estimate for sums of independent Bernoulli variables to obtain a logarithmic bound for $X$ with high probability. By Chernoff’s inequality,
$$\Pr(X\ge t)\le\left(\frac{e\E[X]}{t}\right)^t,$$
see, e.g., Alon–Spencer \cite[Section A.1]{AS}. Take $t=\frac{C\log m}{\log\log m}$ with $C>0$ large. Then
$$\Pr\left(X>\frac{C\log m}{\log\log m}\right)\le\left(\frac{eb}{\frac{C\log m}{\log\log m}}\right)^{\frac{C\log m}{\log\log m}}\le\exp\left(-bC\log m\right).$$
To extend from a fixed great circle to all great circles, choose a collection $\{G_k\}_{k=1}^K$ such that every great circle $G$ is contained in a $m^{-1}$-neighborhood of some $G_k$. One can take $K=O(m^2)$ by selecting the poles of $G_k$ as a maximal $\frac13m^{-1}$-separated set. Consequently, for each $G\subset\Sb$,
$$\Nc_{m^{-1}}(G)\subset\Nc_{2m^{-1}}\left(G_k\right)\quad\text{for some }k=1,...,K.$$
Hence,
\begin{eqnarray*} 
&&\Pr\left(\text{There is a great circle }G\text{ such that }\#\left\{j:p_j\in\Nc_{m^{-1}}(G)\right\}>\frac{C\log m}{\log\log m}\right)\\ 
&\le&\Pr\left(\text{There is a great circle }G_k\text{ such that }\#\left\{j:p_j\in\Nc_{2m^{-1}}\left(G_k\right)\right\}>\frac{C\log m}{\log\log m}\right)\\ 
&\le&K\cdot\exp\left(-bC\log m\right)\\ 
&\le&O\left(m^2\right)\cdot\exp\left(-bC\log m\right)\\ 
&=&o(1)\quad\text{as }m\to\infty, 
\end{eqnarray*}
for some absolute constant $C>0$. This shows that (G) holds for $\{p_j\}_{j=1}^m$ with probability $1-o(1)$ as $m\to\infty$.
\end{proof}

\end{document}